\documentclass{article}

\usepackage{amsfonts}
\usepackage{amssymb}
\usepackage{amsbsy}
\usepackage{amsmath}
\usepackage{amsthm}
\usepackage{url}

\usepackage{citesort}

\newtheorem{theorem}{Theorem}

\newtheorem{definition}[theorem]{Definition}

\newcommand{\Z}{\mathbf{Z}}

\newcommand{\abcd}{\left(\begin{smallmatrix}a & b \\ c &
d\end{smallmatrix}\right)}

\newcommand{\Magma}{{\sc Magma}}
\newcommand{\Sage}{{\sc Sage}}

\DeclareMathOperator{\SL}{SL}

\newcommand{\sltz}{\SL_2(\Z)}

\title{Experimental finding of modular forms for noncongruence subgroups}

\author{L. J. P. Kilford}

\begin{document}

\maketitle

\begin{abstract}

In this paper we will use experimental and computational methods to find
modular forms for non-congruence subgroups, and the modular forms for
congruence subgroups that they are associated with via the
Atkin--Swinnerton-Dyer correspondence. We also prove a generalization
of a criterion due to Ligozat for an eta-quotient to be a modular
function.\end{abstract}

\section{Introduction}
Let~$N$ be a positive integer. We define~$\Gamma(N)$ to be the group
of invertible $2 \times 2$ matrices with coefficients in~$\Z$ whose
reduction modulo~$N$ is congruent to the identity matrix. We say that
a subgroup of the ``modular group'' $\sltz$ is a \emph{congruence
subgroup} if it contains~$\Gamma(N)$ for any~$N$. It can be shown that
these subgroups have finite index in~$\sltz$. We define a
\emph{non-congruence subgroup} to be a subgroup of finite index
inside~$\sltz$ which is not a congruence subgroup.

The theory of modular forms for congruence subgroups is
well-established, at least in integral weights; there are algorithms
to compute bases of space of modular forms, and a well-understood
arithmetic theory of Hecke operators acting on these spaces. There are
many good introductions to this; see~\cite{serre-cours}
or~\cite{diamond-shurman}, for example.

However, as~\cite{atkin-li-long} says in its introduction, the theory of modular forms for non-congruence subgroups is much less well-known, despite the fact that (in a sense that can be
made rigorous) most subgroups of the modular group of finite index are
not congruence subgroups (see~\cite{jones} for a more precise
statement of this result).

The pioneering experimental work of~\cite{atkin-swinnerton-dyer}
discovered congruences satisfied by certain modular forms of this
type, which are now known as Atkin--Swinnerton-Dyer congruences, and
certain of these have been proved by Scholl in a series of
papers~\cite{scholl-deRham,scholl-stdn,scholl-ladic,scholl-ladic2},
which also consider the Hecke algebras attached to spaces of modular
forms for non-congruence subgroups.

More recently, there has been work on refining the conjectured
congruences by Atkin, Li, Long and Yang;
see~\cite{atkin-li-long,li-long-yang,non-congruence-subgroups}. They
prove the Atkin--Swinnerton-Dyer congruences for certain specific
cases, and give another version of the conjectures made earlier in the
field. They also show that the L-functions attached to certain
non-congruence modular forms by Scholl are ``modular'', in the sense
that they can be attached to modular forms for congruence subgroups.

Experimentally, it has been noted by Atkin and Swinnerton-Dyer and
others that the denominators of modular forms for non-congruence
subgroups are unbounded; this is in sharp contrast to the situation
for modular forms for congruence subgroups, which are well-known to
have bounded denominators. It is an interesting open question whether
all  modular forms for non-congruence subgroups have this unbounded
denominator property; recently~\cite{kurth-long} considered this
question; they prove the unbounded denominator property for certain
non-congruence subgroups.

There has also been recent computational work by Verrill et
al~\cite{verrill-et-al}, who have found a number of new examples of
modular forms for non-congruence subgroups which are conjectured to
satisfy the Atkin--Swinnerton-Dyer congruences. This paper inspired
the current work, which also gives lists of modular forms for
non-congruence subgroups which are conjectured to satisfy
Atkin--Swinnerton-Dyer congruences.

The computational work referred to above found congruences involving
modular forms for non-congruence subgroups of genus~0. In
Richards~\cite{richards-thesis}, algorithms are given which extend
this to general non-congruence subgroups, and explicit examples for
genus~1 groups are exhibited. This work is particularly interesting
because it uses complex approximations to modular forms rather than
$p$-adic approximations, thus giving a different and unusual
perspective on the subject.

\section{Notation}

We first give an explicit description of the Atkin--Swinnerton-Dyer
congruence relation, following that given in~\cite{li-long-yang}. We
suppose that~$\Gamma$ is a non-congruence subgroup of finite index in
the modular group~$\SL_2(\Z)$, and that~$k$ is a non-negative integer.

\begin{definition}
\label{asd-relation}
Suppose that~$\Gamma$ has cusp width~$\mu$ at infinity, and that $h
\in S_k(\Gamma)$ has an $M$-integral Fourier expansion at infinity in
terms of~$q^{1/\mu}$ of the form
\[
h(q)=\sum_{n=1}^\infty a_n q^{n/\mu},
\]
for some integer~$M$.

Let~$f$ be a normalized newform of weight~$k$, level~$N$ and
character~$\chi$ (for some congruence subgroup) with Fourier expansion
at infinity given by
\[
f(q)=\sum_{n=1}^\infty c_n q^n.
\]
We say that the forms~$f$ and~$h$ \emph{satisfy the
Atkin--Swinnerton-Dyer congruence relation} if, for all primes~$p$ not
dividing~$MN$ and for all positive integers~$n$, we have that
\begin{equation}
\label{asd-equation}
\frac{a_{np}-c_p a_n +\chi(p)p^{k-1}a_{n/p}}{(np)^{k-1}}
\end{equation}
is integral at all places dividing~$p$. We define~$a_{n/p}$ to be zero
if~$p \nmid n$.
\end{definition}
This is modelled upon the following well-known recurrence relation
that holds for the Fourier coefficients of modular forms for
congruence subgroups which are normalized simultaneous eigenvectors
for the Hecke operators:
\[
a_{np}-a_p a_n +\chi(p) p^{k-1} a_{n/p} = 0,
\]
where again we take~$a_{n/p}$ to be zero if~$p \nmid n$.

Again inspired by the existence of a basis of normalized eigenforms
for spaces of newforms in the congruence case, we now define an
Atkin--Swinnerton-Dyer basis.
\begin{definition}
Let~$k$ be a non-negative integer and let~$\Gamma$ be a non-congruence
subgroup. We say that~\emph{$S_k(\Gamma)$ has an
Atkin--Swinnerton-Dyer basis} if for every prime number~$p$ there is a
basis~$\{h_1,\ldots,h_n\}$ of~$S_k(\Gamma)$ and normalized
newforms~$f_1,\ldots,f_n$ such that each pair~$(h_i,f_i)$ satisfies
the Atkin--Swinnerton-Dyer congruence relation given
in~\eqref{asd-equation}.
\end{definition}
We note that there are cases where one choice of
Atkin--Swinnerton-Dyer basis will suffice for all but finitely many
primes~$p$, and others where the basis depends on the value of~$p$
modulo some integer~$N$. We will describe these below.

\section{Extending the Ligozat criterion}
First, we recall the definition of the Dedekind $\eta$ function;
let~$z$ be an element of the Poincar\'e upper half plane. Then we have
\[
\eta(z)=q^{1/24}\prod_{n=1}^\infty (1-q^n), \text{ where }q:=\exp(2\pi i 
z).
\]
The $\eta$-function can be used to build many interesting modular
forms; for instance, the $\Delta$-function is the $24^{\rm th}$ power
of~$\eta$, and in~\cite{kilford-eta} it is proved that every modular
form for certain congruence subgroups can be written as a sum of
$\eta$-quotients.

We will prove a generalization of the criterion of Ligozat given in
Section~3 of \cite{ligozat} for an $\eta$-quotient to be a modular
function with character for~$\Gamma_0(N)$.

We recall that, if $\abcd\in \SL_2(\Z)$, with~$c \ge 0$, then
\begin{equation}
\label{eq:eta-transform}
\eta\left(\frac{az+b}{cz+d}\right)=\varepsilon\left(\begin{matrix}a &
b\\c & d\end{matrix}\right)\cdot (-i(cz+d))^{1/2} \cdot \eta(z),
\end{equation}
where
\[
\varepsilon\left(\begin{matrix}a & b\\c & d\end{matrix}\right) =
\exp\left(-i\pi\alpha\left(\begin{matrix}a & b\\c &
d\end{matrix}\right)\right)\text{ and }\alpha\left(\begin{matrix}a &
b\\c & d\end{matrix}\right) \in \Z.
\]
The actual definition of~$\alpha$ is rather complicated (it involves
Dedekind sums; for more details, see Section~2.8 of \cite{iwaniec} for
the full story, for instance), but if we have~$(a,6)=1$, then the
following congruence holds:
\[
\alpha\left(\begin{matrix}a & b\\c & d\end{matrix}\right) \equiv
\frac{1}{12} \cdot
a(c-b-3)-\frac{1}{2}\left(1-\left(\frac{c}{a}\right)\right) \mod 2.
\]
This is very useful because it can be shown that~$\Gamma_0(N)$ can be
generated by matrices of the form
\[
\left(\begin{matrix}a & b\\Nc & d\end{matrix}\right)\in\Gamma_0(N),
\text{ with }(a,6)=1\text{ and }a,c \ge 0,
\]
so we need only verify the transformation condition on matrices of
this form to prove our theorem.

Let~$N$ be a positive integer and
define~$g(z)=\prod_{\delta|N}\eta(\delta z)^{r_\delta}$.

\begin{theorem}
\label{ligozat-theorem}
Let~$N$ be a positive integer and let~$g(z)$ be as defined above. Suppose 
that:
\begin{enumerate}

\item\label{x} $\sum_{\delta |N} r_\delta \cdot \delta \equiv 0 \mod 24$,

\item\label{xx} $\sum_{\delta|N} r_\delta \cdot (N/\delta) \equiv 0 \mod 
24$ and

\item\label{xxx} $\sum_{\delta |N} = 0$.

\end{enumerate}
Then~$g(z)$ is a modular function of weight~0 for~$\Gamma_0(N)$ with
quadratic character~$\chi:=\prod_{\delta|N}
\left(\frac{N/\delta}{\cdot}\right)^{r_\delta}$. 
\end{theorem}
\begin{proof}
We take~$U=\abcd\in \Gamma_0(N)$ and~$\delta$ to be a divisor of~$N$.
By explicit computation, we see that
\[
\eta(\delta Uz)=\eta(U_\delta \cdot \delta z)\text{ where }U_\delta=
\left(\begin{matrix}a & b\delta\\c\delta^\prime & d\end{matrix}\right)
\text{ with }\delta \cdot \delta^\prime = N.
\]
Using the explicit formula for the transformation of~$\eta$ given
in~\eqref{eq:eta-transform}, we see that we have
\[
g(Uz)=(-i(Ncz+d))^{\sum_{\delta|N}\frac{r_\delta}{2}} \cdot g(z) \cdot
\prod_{\delta|N} \varepsilon(U_\delta)^{r_\delta}.
\]
From assumption~\eqref{xxx} of the theorem, we see that the first
factor vanishes, so we now need to evaluate the third factor. 
of the cases that we are considering in the theorem we will
have~$(a,6)=1$ (either we have a generator of~$\Gamma_0(N)$ in which
case we can assume this, or we have the auxiliary level
structure~$\Gamma(3)$ which will also allow us to assume this), so we
can rewrite the third factor as
\[
\prod_{\delta|N} \varepsilon (U_\delta)^{r_\delta} =
\exp(-i\pi\lambda)\text{ where }\lambda = \sum_{\delta|N} r_\delta
\cdot \alpha(U_\delta).
\]
Now using the fact that~$(a,6)=1$ (because we are dealing with  a
generator of~$\Gamma_0(N)$) we can write out~$\alpha(U_\delta)$
explicitly as
\[
\alpha(U_\delta) \equiv
\frac{1}{12}a(c\delta^\prime-b\delta-3)-\frac{1}{2}\left(1-\left(\frac{c\delta^\prime}{a}\right)\right)\mod
2,
\]
which means that we can write~$\lambda$ modulo~2 as
\[
\lambda \equiv \frac{1}{12}\left(\sum_{\delta|N} r_\delta \cdot
\delta^\prime\right)-\frac{1}{12}ab\left(\sum_{\delta|N}r_\delta \cdot
\delta\right)-\frac{a}{4}\sum_{\delta|N}r_\delta
-\frac{1}{2}\sum_{\delta|N}\left[1-\left(\frac{c\delta^\prime}{a}\right)\right]\cdot
r_\delta.
\]
We now use the fact that the sum of the degrees~$r_\delta$ is~0 to
show that the third term of the right-hand side is~0. As the
congruences in~\eqref{x} and~\eqref{xx} hold modulo~24, the first and
second terms will vanish modulo~2. 
that our matrix~$\abcd \in \Gamma(3)$ and in particular that~$b\equiv
c\equiv 0 \mod 3$ to show that the first and second terms in the
congruence for~$\lambda$ are integral and still vanish modulo~2.

This means that~$\lambda$ in fact satisfies the congruence
\[
\lambda \equiv 
\frac{1}{2}\sum_{\delta|N}\left[1-\left(\frac{c\delta^\prime}{a}\right)\right]\cdot
r_\delta\mod 2,
\]
and therefore that we have
\[
\exp(-i\pi\lambda) = \prod_{\delta |N} \left(\frac{\delta^\prime
c}{a}\right)^{r_\delta} =
\prod_{\delta|N}\left(\frac{\delta^\prime}{a}\right)^{r_\delta};
\]
where we can take out the factor of~$c$ using~\eqref{xxx}. This means
that we can rewrite~$g(Uz)$ as
\[
g(Uz)=\prod_{\delta|N}\left(\frac{\delta^\prime}{a}\right)^{r_\delta} 
g(z),
\]
so we have shown that~$g$ transforms correctly under the action of
elements of~$\Gamma_(N)$, which proves our theorem.
\end{proof}
The proof will also go through if we take the congruences in \eqref{x}
and~\eqref{xx} modulo~8; in that case, the proof shows that~$g$ is a
modular function for the congruence subgroup~$\Gamma_0(N) \cap
\Gamma(6)$. We cannot expect to prove that it is a modular function
for~$\Gamma_0(N)$ because its Fourier expansion is given in terms
of~$q^{1/3}$ and not~$q$.

We note also that if all of the~$r_\delta$ are even, then the fact
that~$\eta(q)^2$ generates the space of modular
forms~$S_1(\Gamma(12))$ implies that~$g$ is a modular function for the
congruence subgroup~$\Gamma(12N)$.

\section{Algorithm used for finding  modular forms for non-congruence 
subgroups}
The method that we use here is basically a converse to that described
in~\cite{verrill-et-al}. We consider roots of $\eta$-quotients of the
form
\begin{equation}
\label{eta-quotient}
\sqrt[3]{\eta(q^a)^m\eta(q^b)^n\eta(q^c)^r\eta(q^d)^s},
\end{equation}
where the~$a,b,c,d$ are positive integers which divide either~6 or~8,
and~$m+n+r+s=18$ (so the modular function given
in~\eqref{eta-quotient} has weight~3). We assume the unbounded
denominator question discussed in the introduction, that modular forms
for non-congruence subgroups have unbounded denominator, to speed up
the calculations. We cannot use our extension of a theorem of Ligozat here, because it deals with $\eta$-quotients rather than their roots.

We now consider the specific situation where~$p$ is a prime, the
weight~$k$ is~3, $n$ is a positive integer not divisible by~$p$, and
we have an Atkin--Swinnerton-Dyer basis (with respect to~$p$) of our
space of modular forms for a noncongruence subgroup which is composed
of $\eta$-quotients (this is called ``Case 1''
in~\cite{verrill-et-al}). In this particular case, the
equation~\eqref{asd-equation} reduces to
\[
\frac{a_{np}-c_p a_n}{(np)^2},
\]
so in particular we see that~$a_{pn} \equiv c_p a_n \mod p^2$, and
(if~$a_n\ne 0$) then we have~$a_{np}/a_n \equiv c_p \mod p^2$. If we
now fix~$p$ and let~$n$ vary, then the term on the right hand side of
our equation will remain constant, as it does not depend on~$n$, so we
have
\begin{equation}
\label{ap-cp-mod-p2}
a_p \equiv c_p \mod p,
\end{equation}
as long as all of the terms that we have been manipulating were
nonzero modulo~$p^2$. It may happen that we have to consider a
twist~$f \otimes \chi$ of~$f$ to actually get congruences for every
prime; this we can detect by checking to see if~$c_p/a_p$ is a root of
unity for all of the primes~$p$.

However, for some primes the Atkin--Swinnerton-Dyer basis will not be
a pair of distinct $\eta$-quotients~$h_1$ and~$h_2$, but will instead
be of the form~$h_1+\alpha h_2$, where~$\alpha$ is an algebraic number
of small degree (this is called ``Case 2'' in~\cite{verrill-et-al}).
We assume that this means that we have
\begin{equation}
a_{pn}+\alpha b_{pn} \equiv c_p (a_n + \alpha b_n) \mod p^2,
\end{equation}
where as above the~$p^2$ comes from~\eqref{asd-equation} with~$k=3$.
This will hold if
\begin{equation}
\label{ap-constant-mod-p2}
a_{pn} \equiv c_p \alpha b_n \mod p^2 \text{ and }\alpha b_{pn} \equiv
c_p a_n \mod p^2,
\end{equation}
and if every term here is nonzero modulo~$p^2$ then this implies
that~$a_{pn}/b_n \equiv c_p \alpha$ modulo~$p^2$ and~$b_{pn}/a_n
\equiv c_p /\alpha$ modulo~$p^2$. As above, we see that the right-hand
side of these congruences do not depend on~$n$, so if we fix~$p$ and
vary~$n$ then we will find that both~$a_{pn}/b_n$ and~$b_{pn}/a_n$ are
constant modulo~$p^2$.

We can use these congruences to find~$\alpha^2$ and~$c_p^2$
modulo~$p^2$ by combining the congruences above; we find (assuming
that the terms are nonzero) that
\begin{equation}
\label{ap-alpha-mod-p2}
\alpha^2 \equiv \frac{\frac{a_{np}}{b_n}}{\frac{b_{np}}{a_n}} \text{
and } c_p^2 \equiv \frac{a_{np}}{b_n}\cdot\frac{b_{np}}{a_n};
\end{equation}
both of these quantities are (at least experimentally) well-defined
because we have shown that the terms are constant modulo~$p^2$. As in
the previous case, we may need to consider a twist of the form~$f$ by
a character~$\chi$.

We will use this in the following way; we will run over all
$\eta$-quotients of the form~\eqref{eta-quotient} where~$m,n,r,s$ are
less than some bound, and find pairs of $\eta$-quotients which satisfy
one of the the two cases described above for each prime~$p$ up to a
specified bound. The calculations here were performed using
\Magma{}~\cite{magma}; other computer algebra packages such as
\Sage{}~\cite{sage} would also be suitable for this.

Given an experimentally found $\eta$-quotient, we would like to show
that this is a modular form. We will do this by writing it as
the product of a known modular form for a non-congruence subgroup and
a modular function of weight~0 for a congruence subgroup, which will
show that it is a modular \emph{function} for the intersection of
these groups, and then we will verify that its cube is a modular form,
so it has no poles on the upper half plane and therefore is a modular \emph{form}.

The final part of the puzzle is to identify the modular form~$f$ for a
congruence subgroup which satisfies an Atkin--Swinnerton-Dyer
congruence with our $\eta$-quotient. This is mostly a matter of trial
and error; one can guess that the level will be divisible by the
primes~2 and~3, and using \Magma{} we can compute spaces of modular
forms of weight~3 for congruence subgroups. We also have some idea of
what the coefficients of~$f$ should be, because we can
use~\eqref{ap-cp-mod-p2} and~\eqref{ap-alpha-mod-p2} to work out what
those coefficients or their squares are modulo~$p^2$.

\section{Tables of results}
We list some modular forms, mostly taken from Table~12 of~\cite{verrill-et-al} for~$\Gamma_1(6)$ and~$\Gamma_1(12)$, which we can use as building blocks for our non-congruence modular forms.
\begin{eqnarray*}
a&=&\frac{\eta(q)\eta(q^6)^6}{\eta(q^2)^2\eta(q^3)^3}=q-q^2+q^3+q^4+\cdots\\
b&=&\frac{\eta(q^2)\eta(q^3)^6}{\eta(q)^2\eta(q^6)^3}=1+2q+4q^2+2q^3+\cdots\\
c&=&\frac{\eta(q^2)^6\eta(q^3)}{\eta(q)^3\eta(q^6)^2}=1+3q+3q^2+3q^3+\cdots\\
d&=&\frac{\eta(q)^6\eta(q^6)}{\eta(q^2)^3\eta(q^3)^2}=1-6q+12q^2-6q^3+\cdots\\
e&=&\frac{\eta(q)^2\eta(q^3)^2}{\eta(q^2)\eta(q^6)}=1 - 2q - 2q^3+\cdots\\
\end{eqnarray*}

We follow the notation of Verrill et al for these forms; all of them
apart from~$e$ are modular forms of weight~1 for~$\Gamma_1(6)$,
and~$e$ is a modular form of weight~1 for~$\Gamma_1(12)$.

Similarly, there are modular forms and functions listed in Table~11 of~\cite{verrill-et-al} for $\Gamma_1(4) \cap \Gamma_0(8)$ and~$\Gamma_1(16)$ which we can use to construct non-congruence modular forms. Again, we follow the notation given in~\cite{verrill-et-al}.
\begin{eqnarray*}
t&=&\frac{\eta(q)^8\eta(q^4)^4}{\eta(q^2)^{12}} \in M_0(\Gamma_1(4) \cap \Gamma_0(8))\\
\frac{t+1}{2}&=&\frac{\eta(q)^4\eta(q^4)^{14}}{\eta(q^2)^{14}\eta(q^8)^4} \in  M_0(\Gamma_1(4) \cap \Gamma_0(8))\\
\frac{t+1}{2t}&=&\frac{\eta(q^4)^{10}}{\eta(q)^4\eta(q^2)^{2}\eta(q^8)^4} \in  M_0(\Gamma_1(4) \cap \Gamma_0(8))\\
\frac{4(t+1)}{1-t}&=&\frac{\eta(q^4)^{12}}{\eta(q^2)^{4}\eta(q^8)^8} \in  M_0(\Gamma_1(4) \cap \Gamma_0(8))\\
\sqrt{t}&=&\frac{\eta(q)^4\eta(q^4)^2}{\eta(q^2)^6} \in M_0(\Gamma_1(16))\\
\sqrt{\frac{t+1}{2}}&=&\frac{\eta(q)^2\eta(q^4)^{7}}{\eta(q^2)^{7}\eta(q^8)^2} \in  M_0(\Gamma_1(16))\\
E_a&=&\frac{\eta(q^2)^6\eta(q^4)^4}{\eta(q)^4} \in M_3(\Gamma_1(4) \cap \Gamma_0(8))\\
E_b&=& \left(\frac{2t}{t+1}\right) E_a = \frac{\eta(q^2)^8\eta(q^8)^4}{\eta(q^4)^6} \in M_3(\Gamma_1(4) \cap \Gamma_0(8)).
\end{eqnarray*}

Firstly, we present two tables of forms listed in Tables 13 and~14
of~\cite{verrill-et-al}, which have been shown to be modular forms for
certain explicit non-congruence subgroups contained
within~$\Gamma_1(4) \cap \Gamma_0(8)$ and~$\Gamma_1(6)$.
\begin{figure}[h!]
\begin{center}
\begin{tabular}{|c | c | c | c |}
\hline
$h_1$ & & $h_2$ &\\
\hline
 $[ 4, 7,
-4, 11 ]$  & $ \left(\sqrt[3]{b/d}\right) acd$ & $[ -4, 11, 4, 7 ]$ & $ \left(\sqrt[3]{b/d}\right)^2 acd$ \\
$[ 13, -2, -7, 14 ]$ &$ \left(\sqrt[3]{b/c}\right) acd$ & $[ 14, -7,
-2, 13 ]$ &$ \left(\sqrt[3]{b/c}\right)^2 acd$\\

\hline
\end{tabular}
\caption{Forms for subgroups of~$\Gamma_1(12)$ from Table~14
of~\cite{verrill-et-al}.}\label{table-forms-gamma1-12-V}
\end{center}
\end{figure}

\begin{figure}[h!]
\begin{center}
\begin{tabular}{|c | c | c | c |}
\hline
$h_1$ & & $h_2$ &\\
\hline
$[ -8, 20, 2, 4 ]$ & $\sqrt[3]{\frac{t+1}{2}} E_b$ & $[ -4, 22, -8, 8 ]$ & $\sqrt[3]{\frac{t+1}{2}}^2 E_b$ \\
$[ -4, 6, 16, 0 ]$ & $t^{1/3}E_a$ & $[ 4, -6, 20, 0 ]$ & $t^{2/3} E_a$\\
$[ 4, 10, -4, 8 ]$ &  $\sqrt[3]{\frac{t+1}{2t}} E_b$ & $[ 8, -4, 10, 4 ]$ &  $\sqrt[3]{\frac{t+1}{2t}}^2 E_b$\\
$[0,20,-6,4]$ & $\sqrt[3]{\frac{4(t+1)}{1-t}} E_b$ & $[0,16,6,-4]$ & $\sqrt[3]{\frac{4(t+1)}{1-t}}^2 E_b$\\
\hline
\end{tabular}
\caption{Forms for subgroups of~$\Gamma_0(8) \cap \Gamma_1(4)$ from
Table~13 of~\cite{verrill-et-al}.}\label{table-forms-gamma1-8-V}
\end{center}
\end{figure}
We now present two tables of pairs of noncongruence forms which form
Atkin--Swinnerton-Dyer bases that we have found experimentally. We
represent the $\eta$ quotient given in~\eqref{eta-quotient} by the
tuple~$[a,b,c,d]$; if there is one form given twice,
then~$\{h_1,h_2\}$ forms an Atkin--Swinnerton-Dyer basis for all but
finitely many primes~$p$, whereas if there are two distinct forms
given then they form an Atkin--Swinnerton-Dyer basis of the
form~$\{h_1\pm\alpha h_2\}$ for primes satisfying a congruence
condition. 

We see that the new non-congruence modular forms we have discovered can be written as products of modular forms in the same way as those found by Verrill et al, so we can think of these forms we have found as fitting into the same framework as those in~\cite{verrill-et-al}.

\begin{figure}[h!]
\begin{center}
\begin{tabular}{|c | c | c | c |}
\hline
$h_1$ & &$h_2$ &\\
\hline
$[ -8, 13, 8, 5 ]$ & $ \left(\sqrt[3]{e/b}\right) abc$ & $[ 8, 5, -8,
13 ]$  &$ \left(\sqrt[3]{e/b}\right)^2 ace$\\
\hline
\end{tabular}
\caption{Forms for subgroups of~$\Gamma_1(12)$ which form AS-D
bases}\label{table-forms-gamma1-12}
\end{center}
\end{figure}

\begin{figure}[h!]
\begin{center}
\begin{tabular}{|c | c | c | c |}
\hline
$h_1$ & & $h_2$ & \\
\hline
$[ -2, 23, -13, 10 ]$ & $\left(\frac{t+1}{2t}\right)^{1/6}E_b$ & $[ -10, 19, 7, 2 ]$ & $\left(\frac{t+1}{2t}\right)^{5/6}E_b$ \\
$[ 8, -12, 22, 0 ]$  &  $\sqrt{t} E_a $ & $[ -8, 12, 14, 0 ]$ & $\sqrt{t}^5 E_a $\\
$[ 0, -8, 30, -4 ]$ & $ \left(\frac{4(t+1)}{1-t}\right)^{8/3} E_b$& $[ 0, 8, 6, 4 ]$  & $t^{-2/3} E_b$ \\
$[ 2, 17, -11, 10 ]$ & $\left(\frac{t+1}{2}\right)^{1/6} E_a$ &$[ 10, -11, 17, 2 ]$ & $\left(\frac{t+1}{2}\right)^{5/6} E_a$\\
\hline
\end{tabular}
\caption{Forms for subgroups of~$\Gamma_0(16) \cap \Gamma_1(4)$ which
form AS-D bases}\label{table-forms-gamma1-8}
\end{center}
\end{figure}

To show that the new examples we have found are modular forms we check
that we can write them as products of the form~$ f\cdot g$, where~$f$
is a known weight~3 modular form for a noncongruence subgroup (such as
those from Verrill et al's tables), and~$g$ is a modular function of
weight~0 for a congruence subgroup, and that~$(f \cdot g)^3$ is a
modular form. In Section~\ref{section-worked-example} we give a fully
worked-out example of this.

We also give more detail on the forms related by
Atkin--Swinnerton-Dyer congruences; their Fourier expansions
at~$\infty$ and the characters of the congruence forms associated to
them.

\section{Worked example}
\label{section-worked-example}
We will now illustrate how we find two modular forms for a
non-congruence subgroup~$\Gamma$ contained in~$\Gamma_0(16)$. We first
run a search for $\eta$-quotients of the
form~$\sqrt[3]{\eta(q)^m\eta(q^2)^n\eta(q^4)^r\eta(q^8)^s}$ which
experimentally satisfy the congruence conditions discussed after
\eqref{ap-constant-mod-p2}, and we find the following two examples
which do not appear in the tables of~\cite{verrill-et-al}:
\[
H_1:=\sqrt[3]{\frac{\eta(q^2)^{12}\eta(q^4)^{14}}{\eta(q)^{8}}} \text{
and }H_2:=\sqrt[3]{\frac{\eta(q)^8\eta(q^4)^{22}}{\eta(q^2)^{12}}}.
\]
We now notice that we can write~$H_1$ and~$H_2$ as products of a
modular function for~$\Gamma_0(16)$ and~$h_1$ and~$h_2$ from the
second row of Figure~\ref{table-forms-gamma1-8-V}; we have that
\[
H_1= \frac{\eta(q^2)^6}{\eta(q)^4\eta(q^4)^2}\cdot h_2\text{ and } H_2
=  \frac{\eta(q)^4\eta(q^4)^2}{\eta(q^2)^6}\cdot h_1.
\]
We now show that the $\eta$-quotients here are actually modular
functions for $\Gamma_0(16)$ using Theorem~\ref{ligozat-theorem}.
We verify easily that the $\eta$-quotient given above satisfies this
for~$N=16$ although not for~$N=8$, which shows that~$H_1$ and~$H_2$
are modular functions of weight~3 for a suitable noncongruence
subgroup.
Using \Magma{} and \Sage{} it can be checked that $H_1^3,H_2^3 \in
S_9(\Gamma_1(16))$, which verifies that both~$H_1$ and~$H_2$ are in
fact modular forms. 

Finally, we will experimentally determine a classical modular form
which satisfies an Atkin--Swinnerton-Dyer congruence. Figure~\ref{table-worked-example} is a table
which shows the values of~$a_{np}/a_n$ and~$b_{np}/b_n$ modulo~$p^2$,
where an empty space indicates that these numbers are not constant
modulo~$p^2$.
\begin{figure}[h!]
\begin{center}
\begin{tabular}{|c|c|c|}
\hline
$p$ & $a_{np}/a_n$, $b_{np}/b_n$ & $a_{np}/b_n$, $b_{np}/a_n$\\
\hline
5 & 6 &\\
7 &  & 0\\
11 & & 0\\
13 & & 10  \\
17 & 30 & \\
19 & & 0 \\
23 & & 0 \\
29 & -42 & \\
31 & & 0\\
37 & & -70\\
41 & -18 &\\
43 & & 0\\
47 & & 0\\
\hline
\end{tabular}
\caption{Experimentally computed values of~$a_{np}/a_n$
and~$b_{np}/b_n$ modulo~$p^2$, where we take primes~$p \ge 5$ and
positive integers~$n$ with~$pn \le 500$.}\label{table-worked-example}
\end{center}
\end{figure}

Let~$\tau$ be the nontrivial character modulo~4 and let~$f \in
S_3(\Gamma_0(144),\tau)$ be the unique normalized new eigenform with
Fourier expansion at~$\infty$ beginning~$f(q) =
q+6q^5+10q^{13}+O(q^{17})$; this is the twist of the level~16
$\eta$-product~$\eta(q^4)^6$ by the Legendre character modulo~3 ($f$
was found by noticing that the $\eta$-product almost satisfied the
congruence conditions, and then working out which twist actually
worked). Here is the Fourier expansion at~$\infty$ of~$f$ up
to~$O(q^{50})$:
\[
q + 6q^5 + 10q^{13} + 30q^{17} + 11q^{25} - 42q^{29} - 70q^{37} -
18q^{41} + 49q^{49}+O(q^{50}).
\]
It can be seen that the Fourier coefficients of~$f$ are congruent
modulo~$p^2$ to those given in Figure~\ref{table-worked-example}.

We conjecture that the  Atkin--Swinnerton-Dyer basis is $\{H_1,H_2\}$
if~$p\not\equiv 5 \mod 12$, and that it is $\{H_1+H_2,H_1+(p^2+1)
H_2\}$ when~$p\equiv 5 \mod 12$.
\section{Acknowledgements}

I would like to thank Helena Verrill and Ken McMurdy for helpful
conversations, and the Heilbronn Institute for its support while I was
writing this paper. Some of the computations were performed on William
Stein's computer {\tt sage} and I would like to thank him for the use
of this machine.

\bibliographystyle{plain}

\end{document}